\documentclass[graybox]{svmult}

\usepackage{type1cm}        % activate if the above 3 fonts are
\usepackage{makeidx}         % allows index generation
\usepackage{graphicx}        % standard LaTeX graphics tool
\usepackage{multicol}        % used for the two-column index
\usepackage[bottom]{footmisc}% places footnotes at page bottom

\usepackage{newtxtext}       % 
\usepackage{newtxmath}       % selects Times Roman as basic font

\makeindex             % used for the subject index
\begin{document}
\title*{Stochastic Solutions of Stefan Problems with General Time-Dependent Boundary Conditions}
%\title*{Stochastic solutions of Stefan problems} %% with general time-dependent boundary conditions}
\author{Magnus \"{O}gren}
\institute{Magnus \"{O}gren, \email{magnus.ogren@oru.se} \at School of Science and Technology, \"{O}rebro University, 70182 \"{O}rebro, Sweden. \\
Hellenic Mediterranean University, P. O. Box 1939, GR-71004, Heraklion, Greece.}
%%\institute{Hellenic Mediterranean University, P. O. Box 1939, GR-71004, Heraklion, Greece}

%%\and Name of Second Author \at Name, Address of Institute \email{name@email.address}}
\maketitle

%%\abstract*{Each chapter should be preceded by an abstract (no more than 200 words) that summarizes the content. The abstract will appear \textit{online} at \url{www.SpringerLink.com} and be available with unrestricted access. This allows unregistered users to read the abstract as a teaser for the complete chapter.
%%Please use the 'starred' version of the \texttt{abstract} command for typesetting the text of the online abstracts (cf. source file of this chapter template \texttt{abstract}) and include them with the source files of your manuscript. Use the plain \texttt{abstract} command if the abstract is also to appear in the printed version of the book.}

\abstract{This work deals with the one-dimensional Stefan problem with a general time-dependent boundary condition at the fixed boundary. 
Stochastic solutions are obtained using discrete random walks, and the results are compared with analytic formulae when they exist, otherwise with numerical solutions from a finite difference method. 
The innovative part is to model the moving boundary with a random walk method. 
The results show statistical convergence for many random walkers when $\Delta x \rightarrow 0$.
Stochastic methods are very competitive in large domains in higher dimensions and has the advantages of generality and ease of implementation. 
The stochastic method suffers from that longer execution times are required for increased accuracy.
Since the code is easily adapted for parallel computing, it is possible to speed up the calculations. 
Regarding applications for Stefan problems, they have historically been used to model the dynamics of melting ice, and we give such an example here where the fixed boundary condition follows data from observed day temperatures at \"{O}rebro airport. 
Nowadays, there are a large range of examples of applications, such as climate models, the diffusion of lithium-ions in lithium-ion batteries and modelling steam chambers for petroleum extraction.
}

\section{Introduction}

The Stefan Problem has its name from Josef Stefan (1835-1893) who was first to investigate problems including a moving boundary in detail. This was described in his report on ice formation in polar seas~\cite{stefan}, where he also presented the analytical solution, see Eq.~(\ref{stefanSol}), to the problem where the fixed boundary has constant temperature. However, for the general case where the temperature at the fixed boundary is an arbitrary function, 
%%Stefan could not give an explicit solution except for some ideas and approximations. 
no explicit solution has been obtained, though there are power series formulations described in the literature, see e.g.~\cite{tao}. 
In addition to the ice formation problem originally examined by Stefan, moving boundary problems now have many other applications, see e.g.~\cite{hunke, hoffman, chen_2015}. 

The aim here is to show how to solve Stefan problems for arbitrary boundary conditions using stochastic methods.
We will present a discrete Random Walk Method (RWM) that solves the Stefan problem, which is a PDE consisting of the heat equation defined in a phase changing medium. There are different types of formulations of this problem, but one of the characteristics is that it has a \textit{free} or \textit{moving boundary} governed by a so-called Stefan condition, which describes the position of the interface between the phases. Beyond the moving boundary, the general formulation of the problem usually also includes a fixed boundary with a boundary condition different from the moving one. For physical reasons the boundary condition at the moving boundary is here set to be the transition temperature, i.e. the melting point $T_M$ of ice.
At the fixed boundary, the condition for the temperature may be set to an arbitrary function $f(t)$ of time. 
For the case where we have a constant temperature at the fixed boundary $f(t)=T_0$ and for one other special form of $f(t)$, there are analytical solutions to the Stefan problem. 
In addition a specific time dependent incoming heat flux is illustrated to be equivalent with the constant temperature condition.
However, in most cases we need numerical calculations to evaluate a solution. As a practical example of such a case, we model the melting of ice where the surface temperature is defined according to the variations in the air temperature. %\cite{umea}. \newline

\subsection{Random walk and the heat equation}

In this Section, we will study the heat equation %% physical constants related to the
\begin{equation}
    \frac{\partial T}{\partial t}= \alpha \frac{\partial^2 T}{\partial x^2} \, ,
    \label{eq:heat_equation_no_bc}
\end{equation}
and describe how to translate it into an RWM~\cite{Chati_2001, OgrenEPJB2014}. 
In our one-dimensional model, we want to let one walker represent the temperature difference of 1$^\circ$C on the volume element $\Delta x \cdot 1 \cdot 1 \hspace{5pt} m^3$. 
To make a simple illustration for the heat equation, denote the number of walkers in the volume element $i$ with width $\Delta x$ at time $t$ as $N_i(t)$. If we let the probability for a walker to go either to the left or to the right to be equal during a time step $\Delta t$ we have equal probabilities $P= 1/2$. Then we expect to have $P  N_i$ walkers going to volume element $i+1$ and the same amount going to volume element $i-1$. At the same time walkers from volume elements $i-1$ and $i+1$ will walk into the volume element $i$, giving the following balancing equation for $N_i$
\begin{equation}
\begin{array}{c}
    \vspace{6pt}
     N_i(t+\Delta t) = N_i(t) - (\text{walkers to } i-1) - (\text{walkers to } i+1) \\
     \vspace{6pt}
      + (\text{walkers from } i-1) + (\text{walkers from } i+1) \\
      \vspace{6pt}
      \Rightarrow N_i(t+\Delta t) - N_i(t) = P ( N_{i+1}  - 2 N_i + N_{i-1} )  \, .
\end{array}
\label{eq:derivation_finite_differences}
\end{equation}
We divide Eq.~(\ref{eq:derivation_finite_differences}) by $\Delta t$ and introduce the constant 
\begin{equation}
\alpha= (P/\Delta t) (\Delta x)^2,
\label{eq:alpha}
\end{equation}
such that
\begin{equation}
    \frac{N_i(t+\Delta t)-N_i(t)}{\Delta t} = \alpha \frac{N_{i+1}(t)-2N_i(t)+N_{i-1}(t)}{(\Delta x)^2} \, ,
    \label{dN}
\end{equation}
which is a discretized partial differential equation for $N(x,t)$.
We then see that Eq.~(\ref{dN}) has the same form as the heat equation~(\ref{eq:heat_equation_no_bc}).

In general, the same arguments can be made to derive the corresponding equation in $D$ dimensions, since a symmetric Cartesian grid has $2D$ directions for a walker to go with equal probability $P=1/(2D)$, hence $\alpha= (\Delta x)^2/(2D\Delta t) $. %% ~\cite{OgrenEPJB2014}. 

We now first present an introductory example without boundary conditions. 
Consider the heat conduction problem for an infinite rod, with a central heat impulse at $t=0$
\begin{align}
    \label{ex1}
    \frac{\partial T}{\partial t}&= \alpha\frac{\partial^2 T}{\partial x^2} \, ,        && x \in \mathbb{R} \, ,        && t>0 \, , \\[6pt]
    T(x,0)&=\delta (x) \, . \label{eq:initial_b}
\end{align}
%
%%From~\cite{sparr} we can see that the 
The well known solution to this problem is
\begin{equation}  %%\label{ana1}
    T(x,t) = \frac{1}{\sqrt{4 \pi \alpha t}} e^{-x^2/(4 \alpha t)} \, .
    \label{eq:sol_no_bc}
\end{equation}
\begin{figure}[t]
    \includegraphics[width = 60mm]{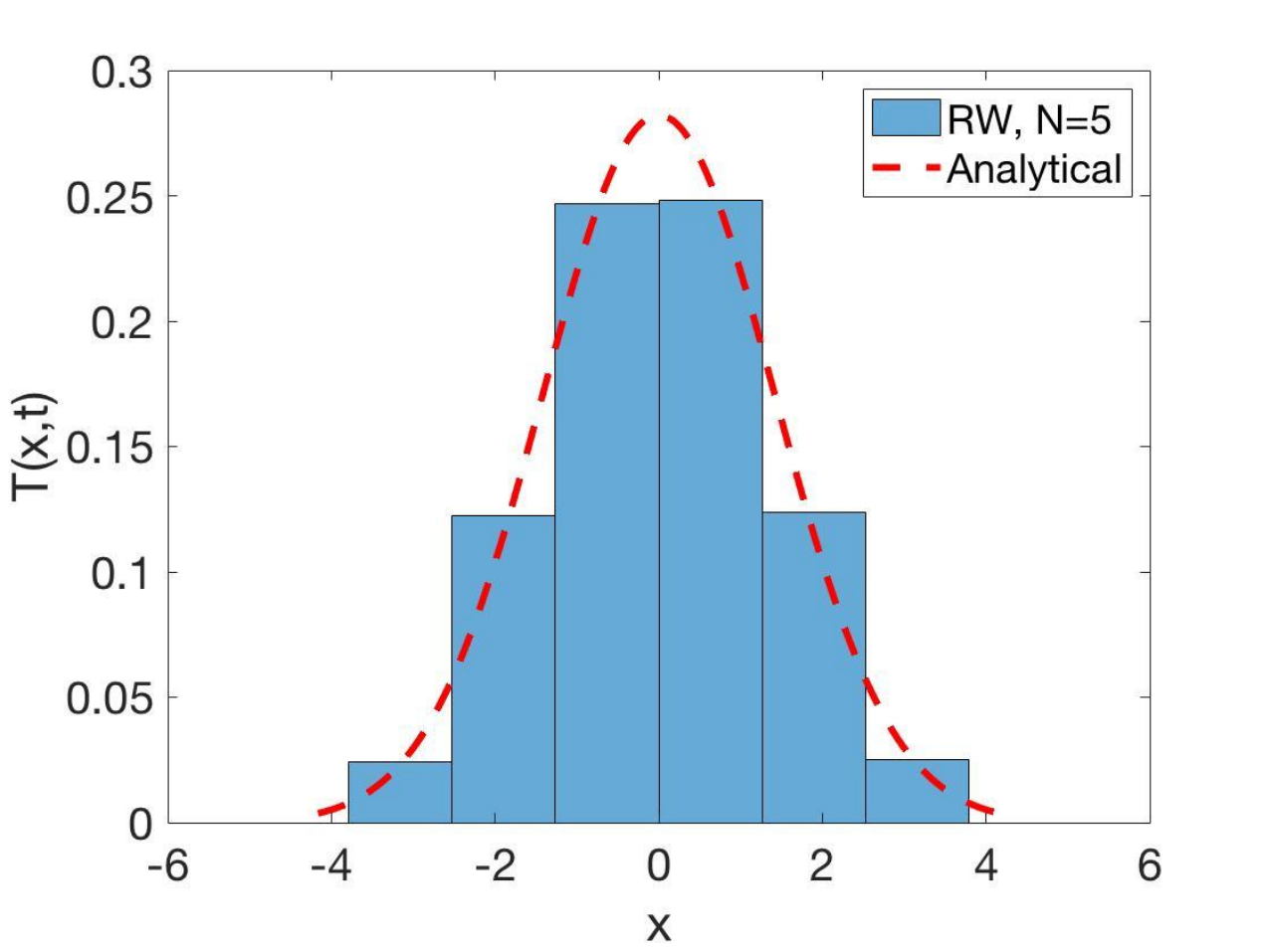}~\includegraphics[width = 60mm]{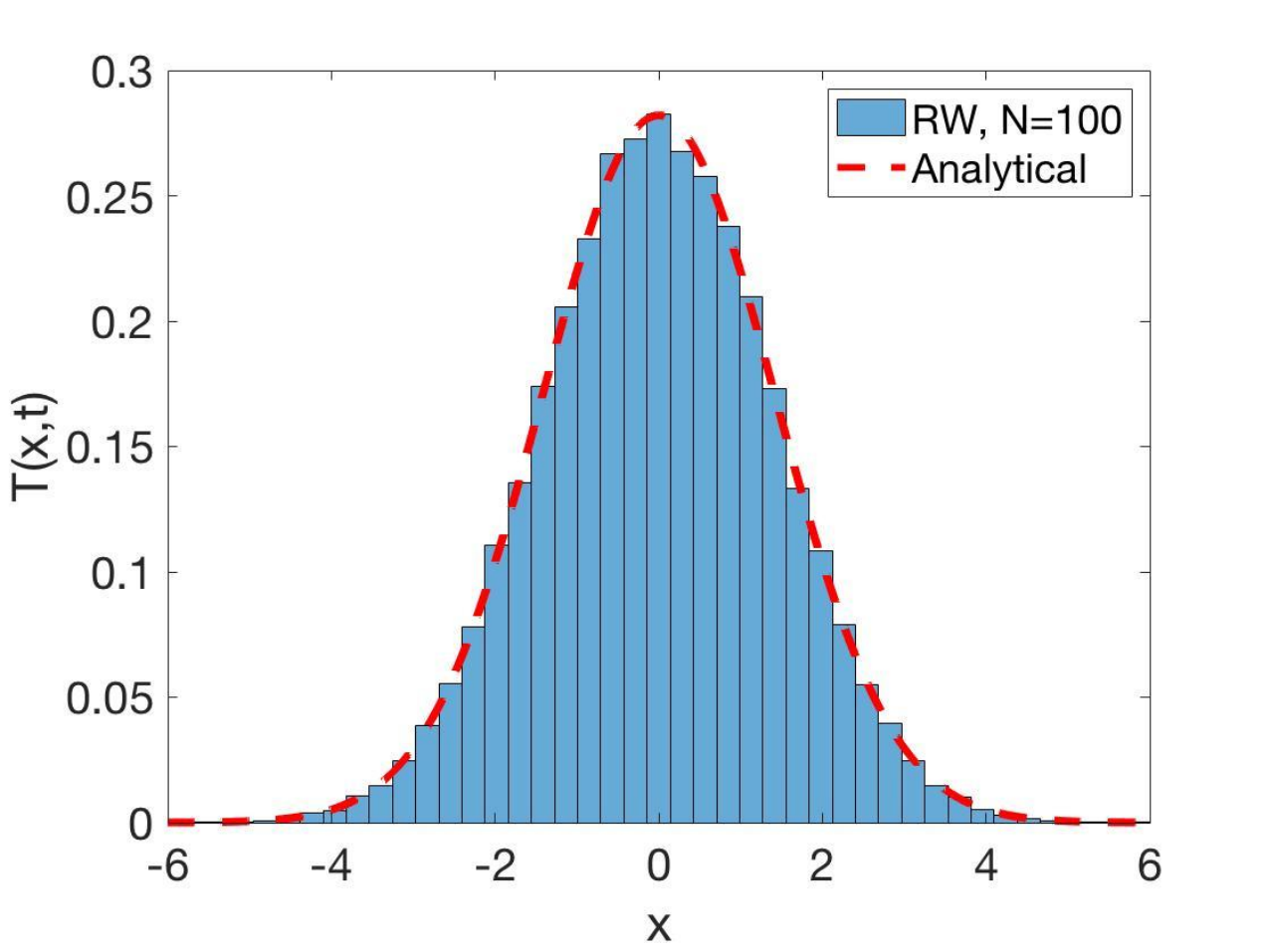}
    \caption{Analytical vs. random walk solutions for the introductory problem of Eqs.~(\ref{ex1})-(\ref{eq:initial_b}) for $t=1$.
    The red dashed curves are from Eq.~(\ref{eq:sol_no_bc}).
    The blue bins have the width $2 \Delta x$.
    (Left) N=5 time steps, i.e. $\Delta t = 1/5$ and $\Delta x = \sqrt{2/5} \simeq 0.63$. 
    (Right) N=100 time steps, i.e. $\Delta t = 1/100$ and $\Delta x = \sqrt{2}/10 \simeq 0.14$.}
    \label{fig_RW}
\end{figure}
This problem is straightforward to model with the RWM since it is defined for all real values and thus has no boundary conditions to consider. 
In Fig. \ref{fig_RW} we can see a comparison between the analytic solution and the discrete probability density function for random walks with $n=10^5$ initial walkers at $x_0=0$.

\subsection{A random walk model with boundary conditions} \label{sec:RW_BC}

As the next problem, a heat equation with fixed boundaries and homogeneous Dirichlet conditions is considered.
\begin{align}
    \label{start1}
    \frac{\partial T}{\partial t}&= \alpha\frac{\partial^2 T}{\partial x^2} \, ,        && 0<x<L \, ,        && t>0 \, , \\[6pt]
    T(0,t)&= T(L,t) = 0 \, ,        && t>0 \, ,       \\[6pt]%\label{start3}
    T(x,0)&= g(x) \, . \label{eq:start_3_b}
\end{align}
Using separation of variables on Eq.~(\ref{start1}), the general solution to this problem can be written %% is~\cite{sparr}
\begin{equation}
    T(x,t) = \sum\limits_{n=1}^{\infty} c_n e^{-\alpha (\frac{n\pi}{L})^2 t} \sin{\left(\frac{n\pi x}{L}\right)} \, .
    \label{eq:_general_sol_fourier}
\end{equation}
We write the initial condition as
\begin{equation}
    g(x) = \sum\limits_{n=1}^{\infty} c_n  \sin{\left(\frac{n\pi x}{L}\right)} \, .
\end{equation}
By recognizing this as the Fourier series expansion of $g(x)$ on $0<x<L$, $c_n$ can be determined according to
\begin{equation}
    c_n = \frac{2}{L} \int\limits_0^L g(x) \sin{\left(\frac{n\pi x}{L}\right)}dx \, .
\end{equation}
Here we set $g(x)=1$ as an example, which gives the solution from Eq.~(\ref{eq:_general_sol_fourier}) on the form
\begin{equation}\label{sol_fourier}
    T(x,t) = \frac{4}{\pi}\sum\limits_{k=0}^{\infty} \frac{\exp{ \left( -\frac{\alpha\pi^2}{L^2}(2k+1)^2 t \right) }}{2k+1} \sin{\left(\frac{(2k+1)\pi x}{L}\right)} \, .
\end{equation}

In the previous problem of Eqs.~(\ref{ex1})-(\ref{eq:initial_b}) we adapted a RWM to a problem defined on the whole x-axis. If we instead want to solve the problem of Eqs.~(\ref{start1})-(\ref{eq:start_3_b}), it is necessary to implement boundary conditions. 
This is done by discretizing the space and time on the finite domain according to Eq.~(\ref{eq:alpha}). 
In this example we choose the following discretization, where we for simplicity set $\alpha=1$
\begin{equation}
    \begin{array}{lll}
        \vspace{6pt}
         \Delta x = x_{i+1}-x_i=10^{-2} \, , & i=0,1, \hdots,N-1 \, , & x_0 = 0 \, , \quad x_N=1 \, , \\
         \vspace{6pt}
         \Delta t = t_{j+1}-t_j=5 \cdot 10^{-5} \, , & j=0,1, \hdots,M-1 \, , & t_0 = 0 \, , \quad t_M=1 \, . \\
    \end{array}
\end{equation}
The initial condition of Eq.~(\ref{eq:start_3_b}) with $T(x,0)=1$ will here be represented by one walker starting at $T(x_i,t_0)$ for all $i$. 
By the next timestep $t_1$, all walkers will have moved one step either to the right or to the left. 
For homogenous Dirichlet conditions, the walkers that reach the boundaries will be absorbed and disappear, such that $T(x_0,t_j) = T(x_N,t_j) = 0$. 
In the case of inhomogenous Dirichlet conditions, i.e. $T(x_0,t_j) \neq 0$, as in the upcoming Stefan problem, see Eq.~(\ref{BCs}), we also have walkers starting from the boundary. The number of walkers starting at $T(x_0,t_j)$ will here be set according to $f(t)$, where $f(t)=1$ will be represented by one walker starting at $T(x_0,t_j)$ for all $j$. 
We then iterate over time until all walkers have reached the boundaries or the maximum time $t_M$ is attained. 
A statistical problem so far is that the result of our model with \textit{one} walker, representing a temperature difference of $1^\circ$C per volume unit, might differ a lot depending on how each random walk turns out. 
Real moving particles causing thermal diffusion representing that raise of temperature are large in numbers. 
Therefore, to get an accurate result we multiply the number of walkers starting at all points defined by initial- or boundary conditions with a large number $n$, and at the end we divide the temperature at all points with $n$.

\begin{figure}[t]
    \includegraphics[width = 60mm]{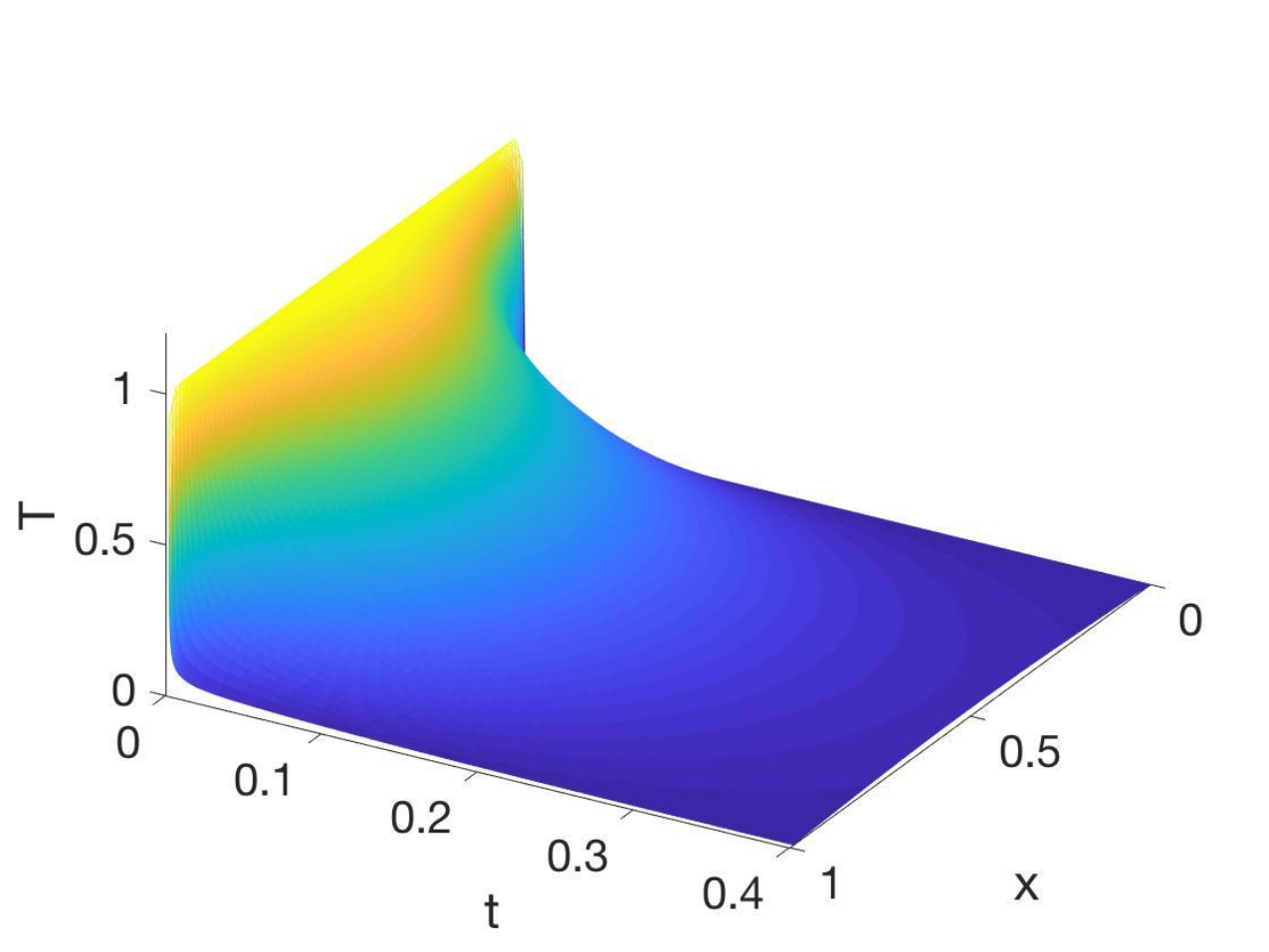}~\includegraphics[width = 60mm]{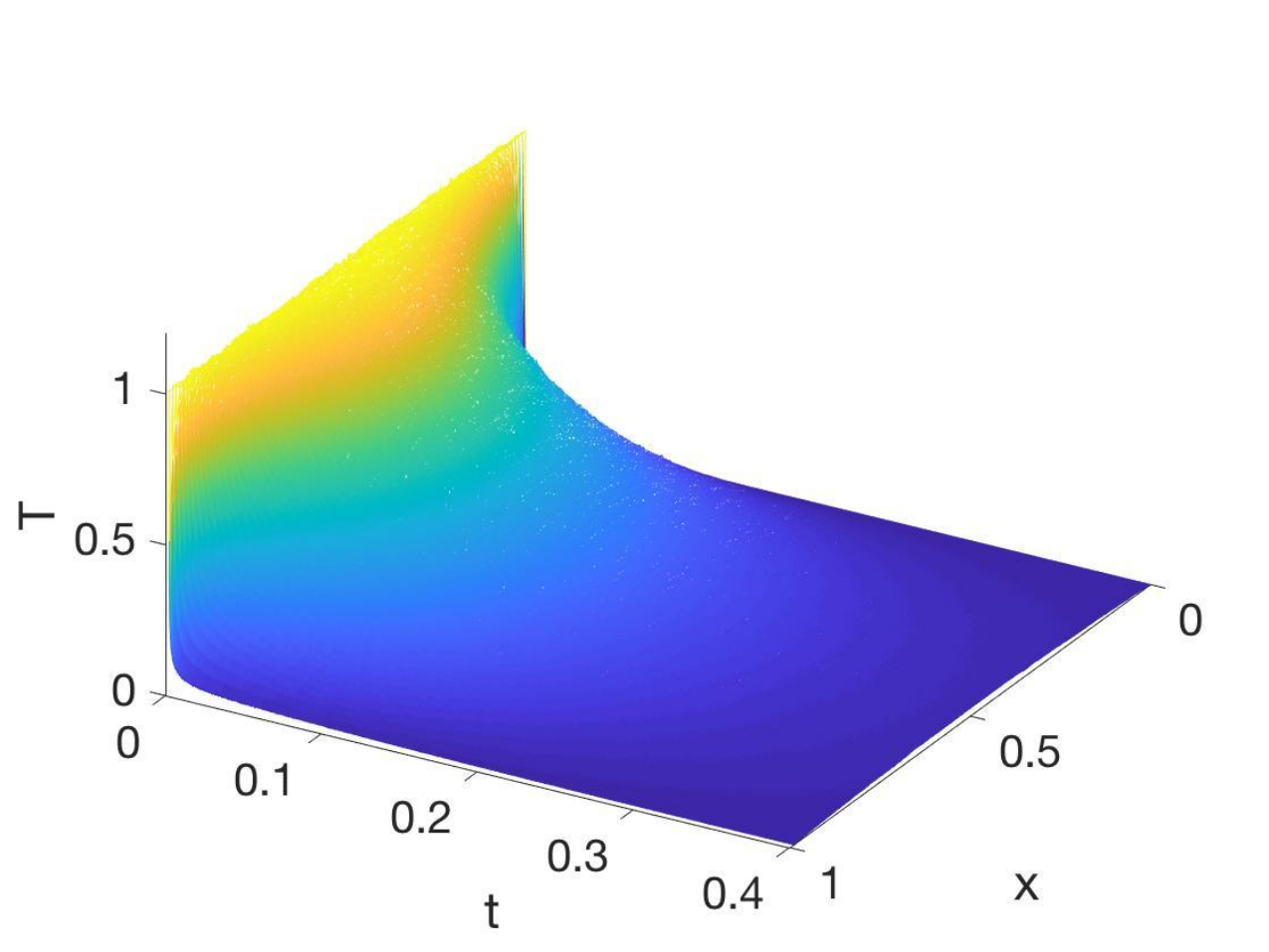}
     \caption{Solutions for the heat conduction problem for $t \in [0,0.4]$ with initial condition $g(x)=1$. 
     (Left) Analytic solution from Eq.~(\ref{sol_fourier}) with 100 terms in the Fourier series. 
     (Right) RWM solution with $n=10^4$ and $\Delta x=0.01$.
     }
    \label{fig_fixed}
\end{figure}
\begin{figure}[t]
\sidecaption[t]
     \includegraphics[width=75mm]{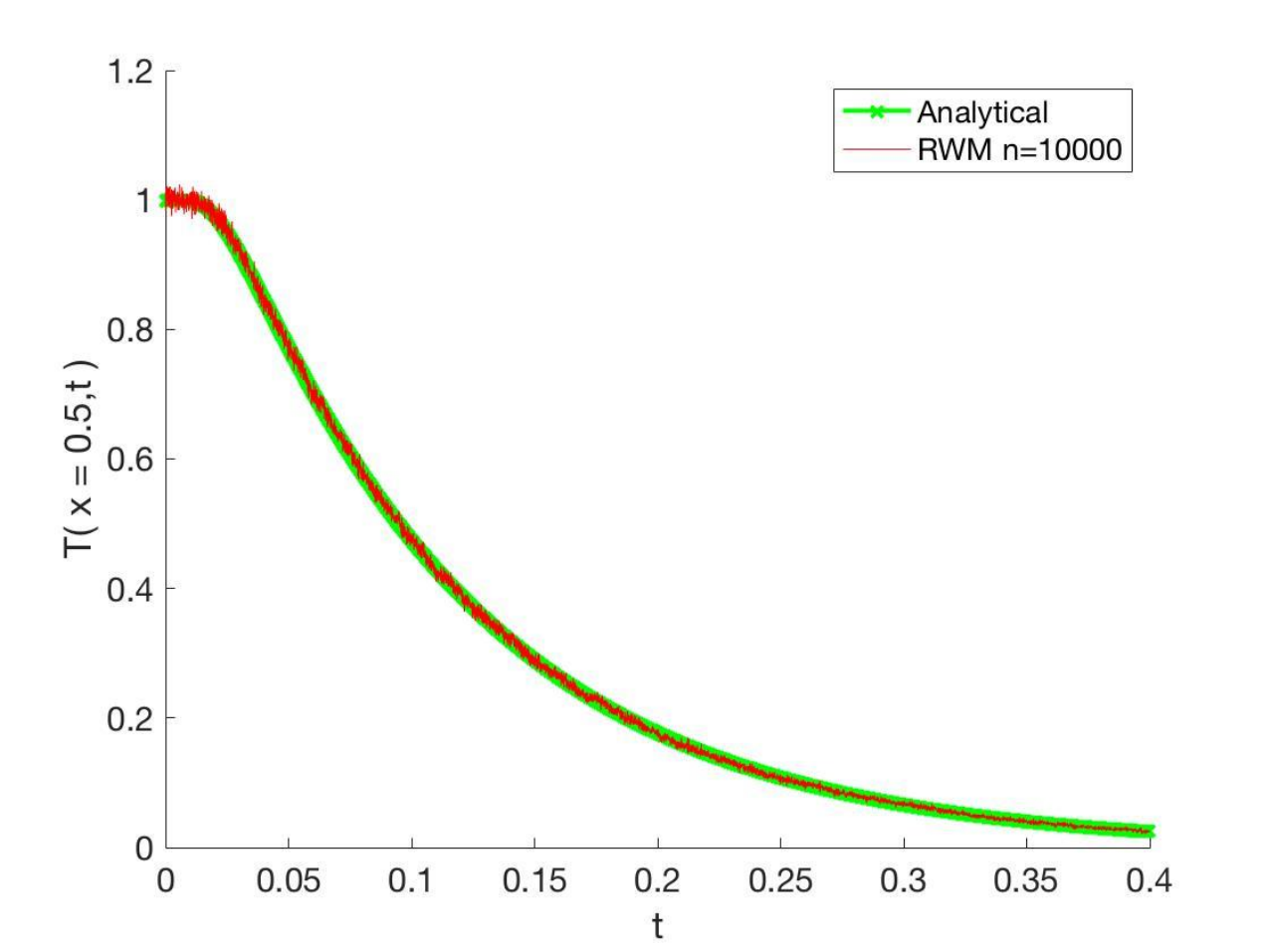}
     \caption{Comparison between the analytic solution, Eq.~(\ref{sol_fourier}) (green) for $x=0.5$ and $t \in [0,0.4]$,  
     and the RWM (red) with $n=10^4$ and $\Delta x=0.01$. }
    \label{fig_fixed_b}
\end{figure}

In Eq.~(\ref{sol_fourier}) we presented an analytic solution for the heat conduction problem Eqs.~(\ref{start1})-(\ref{eq:start_3_b}) for the initial temperature $g(x)=1^\circ$C. 
Fig. \ref{fig_fixed} shows the temperature distributions $T(x,t)$ for the analytical result and the RWM solution with $\alpha=1$. In Fig.~\ref{fig_fixed_b} we see a comparison between the analytical result and the RWM  in the cross-section $x=0.5$.

\section{The Stefan problem}

In our model for the one-dimensional Stefan problem we consider an initial block of ice, i.e. a solid (S), with semi-infinite extent ($x \rightarrow \infty$) and one surface to air at $x=0$. 
At $t=0$ there is no water phase and the temperature for the ice phase is at $T_{ice}= T_m =0^\circ$C. 
For $t>0$ the ice can start to melt and thus we can have a water phase, i.e. a liquid (L), to the left of the ice.  
We presently treat only the so-called one-phase Stefan problem, which means that the temperature in the ice phase does not change in time.
The temperature at the $x=0$ boundary, i.e. the interface between air and water for $t>0$, is allowed to change over time according to $f(t)$, and to simulate a melting process, we initially assume $f(t)>0$. This yields the following equations~\cite{Gupta2003}
\begin{align}
    \label{heat1}
    \frac{\partial T}{\partial t}&= \alpha_L\frac{\partial^2 T}{\partial x^2} \, ,      && 0<x<s(t) \, ,        & t>0 \, , \\[6pt]
    \label{BCs}
    T(0,t)&= f(t) \, ,      && t>0 \, , \\[6pt]
    T(x,0)&=0 \, , \\
    \label{sc}
     \rho l \frac{ds}{dt}&= - k_L \left.\frac{\partial T}{\partial x}\right|_{x=s(t)} \, ,      && t>0 \, , \\[6pt]
     s(0)&=0 \, , \\[6pt]
     \label{ste1}
     T(s(t),t)&= T_M =0 ^\circ \text{C}\, ,        && t>0 \, .
\end{align}
Here the thermal diffusivity in the liquid part, $\alpha_L$ [m$^2$/s] in (\ref{heat1}), is defined as 
\begin{equation}
    \alpha_L = \frac{k_L}{\rho_L c_L} \, ,
\end{equation}
where $k_L$ [W/(mK)] is the heat conductivity, $\rho_L$ [kg/m$^3$] the density and $c_L$ [J/(kgK)] the specific heat capacity in the liquid phase. 
Note that these physical properties differ between the solid and liquid part, e.g. $\alpha_L \neq \alpha_S$.
But since $T=0^\circ$C in the solid phase and the temperature distribution only is evaluated in the liquid phase, $\alpha_S$ is not taken into consideration in this one-phase Stefan problem. 
In equation (\ref{sc}) $l$ is the specific latent heat and $\rho$ is the density. 
Here it is assumed that $\rho = \rho_L = \rho_S$ for simplicity.
The analytic solution to the problem when $f(t)=T_0$ is constant is~\cite{stefan}
\begin{equation}\label{stefanSol}
  %\[
    \left\{
                \begin{array}{ll}
                  \vspace{6pt}
                  T(x,t) & =  T_0\left(1 - \frac{\text{erf}(x/(2 \sqrt{\alpha_L t}))}{\text{erf}(\lambda)}\right) \,  ,  \\
                  \vspace{6pt}
                  s(t) & = 2 \lambda \sqrt{\alpha_L t} \,  , \\
                  \beta \sqrt{\pi} \lambda e^{\lambda^2}\text{erf}(\lambda) & =  T_0 \, ,
                \end{array}
              \right.
  %\]
  \end{equation}
  where $\beta=l/c_L$ and $\text{erf}{(x)}$ is the \textit{error function} defined as $\frac{2}{\sqrt{\pi}} \int_{0}^{x}e^{-y^2} \, dy$. 
  
  A different special case when an analytic solution also exist is when $f(t)=e^t-1$. 
  Provided $\beta=1$, the solution is then~\cite{stefan, vynnycky}
  \begin{equation}\label{special}
          \left\{
                \begin{array}{ll}
                  \vspace{6pt}
                  T(x,t) & =  e^{t-x}-1 \\
                  \vspace{6pt}
                  s(t) & = t \, . \\
                \end{array}
              \right.
  \end{equation}
  We will use also this case for a numerical comparison with the RWM in Sect. \ref{results}.

  \subsection{The Stefan condition}\label{stefan_c}
  
  The position of the free boundary, i.e. the interface between the two phases, is time-dependent and denoted as $x=s(t)$. 
  At time $t_0$ the entire domain $x>0$ is divided into two subdomains consisting of, the water phase $x<s(t_0) $, and the ice phase $x>s(t_0)$. 
  Here we consider a \textit{one-phase problem} which means that the temperature in one of the phases (here the ice phase) is constant at the melting temperature $T_M=0 ^\circ$C. 

\begin{figure}[t]
\sidecaption[t]
  %   \centering
      \includegraphics[width=75mm]{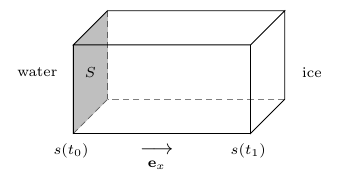}
      \caption{Volume of melted ice between time $t_0$ and $t_1$.}
      \label{fig_melt}
  \end{figure}

  We here briefly derive the Stefan condition stated in Eq.~(\ref{sc}), which will later be used in the formulation of the stochastic model for the interface $s(t)$. 
  More details on the derivation of the Stefan condition can be found e.g. in~\cite{Gupta2003}. 
  In the case of melting ice, the water phase at time $t_1 > t_0$ will be increased, resulting in $s(t_1)>s(t_0)$. 
  If we imagine a block of ice with cross sectional area $S$, the volume $V$ of the melted ice in the time interval $t \in [t_0,t_1]$ is $S  \left( s(t_1)-s(t_0) \right)$, see Fig.~\ref{fig_melt}. 
  The thermal energy $Q$ [J] required for the melting of this block is determined according to
  \begin{equation}\label{Q1}
      Q=V  \rho  l = S \left(s(t_1)-s(t_0) \right) \rho l \, ,
  \end{equation}
  where $l$ [J/kg] is the specific latent heat for the phase transition. 
  As we here assume that the heat is only spread by diffusion, the heat transport obeys Fourier's law
  \begin{equation}
      q=-k_i\frac{dT}{dx} \, ,
  \end{equation}
  where $q$ is the local heat flux density $[\text{W/m}^2]$. 
  By energy conservation and the expressions for the heat fluxes from the liquid and solid phases, $Q$ can be written
  \begin{multline}\label{Q2}
      Q=\int\limits_{t_0}^{t_1} \int \left( -k_L \frac{\partial T(s(\tau),\tau)}{\partial x} \cdot \mathbf{e}_x -k_S \frac{\partial T(s(\tau),\tau)}{\partial x} \cdot (-\mathbf{e}_x) \right) \, dS \, d\tau 
      \\
     = S \int\limits_{t_0}^{t_1}  \left( -k_L \frac{\partial T(s(\tau),\tau)}{\partial x} +k_S \frac{\partial T(s(\tau),\tau)}{\partial x} \right) \, d\tau \, .
  \end{multline}
  Combining Eqs.~(\ref{Q1}) and (\ref{Q2}), dividing by $t_1-t_0$, and letting $t_1 \rightarrow t_0$, will yield Eq.~(\ref{sc}) for the Stefan condition
  \begin{multline} \label{sc1}
      \rho l S \lim_{t_1 \to t_0} \frac{s(t_1)-s(t_0)}{t_1-t_0} = S \lim_{t_1 \to t_0} \frac{1}{t_1-t_0} \int\limits_{t_0}^{t_1}  \left( -k_L \frac{\partial T(s(\tau),\tau)}{\partial x} +k_S \frac{\partial T(s(\tau),\tau)}{\partial x} \right) \, d\tau
       \\
    \Rightarrow \rho l \frac{ds}{dt} = -k_L \frac{\partial T(s(t),t)}{\partial x} +k_S \frac{\partial T(s(t),t)}{\partial x} \, .
  \end{multline}
  Here $t_0$ have been replaced by $t$ since $t_0$ can be chosen arbitrarily. 
  In the present case where we assume $T=0 ^\circ$C for $x>s(t)$, diffusion only occur in the liquid phase and Eq.~(\ref{sc1}) reduces to 
  \begin{equation}
      \rho l \frac{ds}{dt} = -k_L \frac{\partial T(s(t),t)}{\partial x} \, .
  \end{equation}

\subsection{Modelling the moving boundary}

To be able to solve the Stefan problem with the RWM, the critical part is how to handle the moving boundary $s(t)$. To set up a model for the movement of the boundary $s(t)$ we start from Sect. \ref{stefan_c}. 
In Eq.~(\ref{Q1}) we established that the heat required to move the boundary a small step $\Delta s$ is 
\begin{equation} %\label{st1}
Q= S \Delta s  \rho l \, ,
\end{equation}
and thus
\begin{equation}\label{st1}
    \Delta s= \frac{Q}{S \rho l}  \, . 
\end{equation}
\begin{figure}[t]
\sidecaption[t]
      \includegraphics[width=40mm]{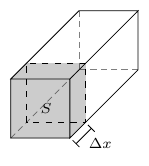}
      \caption{One walker raises the temperature of the gray volume with $1^\circ$C. }
      \label{fig_volym}
\end{figure}
We want to compare this with the heat represented by one walker as it raises the temperature $1 ^\circ$C of the volume $S \Delta x $ [m$^3$], see Fig. \ref{fig_volym}. 
This can be expressed as ($c=c_L$)
\begin{equation}\label{st2}
    Q_{walker} = c  \rho  V  \Delta T = c  \rho S \Delta x   \cdot 1 ^\circ \text{C} \, .
\end{equation}
By combining Eqs. (\ref{st1}) and (\ref{st2}), we have
\begin{equation}
         \Delta s = \frac{c \rho S \Delta x }{l \rho S} = \frac{c}{l} \Delta x \, .
         \label{eq:Delta_s}
\end{equation}
So for every walker absorbed by the moving boundary at $s(t)$ the boundary will move the increment $\Delta s$~\cite{stoor}. 
To adjust for the multiplication with the factor $n$ at the starting points, as discussed in Sect.~\ref{sec:RW_BC}, we also need to correct the step length $\Delta s$ by dividing with $n$. 
Hence, the moving boundary will have the position $i$ in the x-grid when %th esum of all  $s_n = s{n-1} + \Delta s$ is 
\begin{equation}
    i \leq \frac{s_k}{\Delta x} < i+1 \,  , \   s_k = s_{k-1}+\frac{\Delta s}{n} \, .
    \label{eq:Delta_s_k} 
\end{equation}
It is of interest to see how the ratio between $\Delta s$ and $\Delta x$ turns out as we insert realistic physical parameter values for $c$ and $l$. 
For water at $0 ^\circ$C we have $c=4.22$ kJ/(kg$\cdot$K) and $l=334$ kJ/kg~\cite{lockby}, which gives $c/l \approx 0.0126 $,
and we see from Eq.~(\ref{eq:Delta_s}) that $\Delta s \ll \Delta x$. Note that in the opposite case, if $\Delta s \gg \Delta x$, the boundary will move several $\Delta x$-steps as it is reached by one walker and this will lead to poor results when modelling the movement of the boundary. 
Thus, in the case that we have $c/l > 1$ we have to compensate by increasing the number $n$ and thereby decreasing the step size $\Delta s /n$ in Eq.~(\ref{eq:Delta_s_k}). 
So a rule of thumb to yield a good approximation of the boundary is to choose $n$ such that $ \Delta s/ n \ll \Delta x  $.

\subsection{Stefan problem with an incoming heat flux}

In the Stefan problem Eqs.~(\ref{heat1})-(\ref{ste1}) the temperature at the fixed boundary ($x=0$) is described by the Dirichlet condition of Eq.~(\ref{BCs}).
Changing instead to a Neumann condition  
\begin{equation}
    \frac{\partial T}{\partial x} (0,t) = h(t)  \, ,
        \label{eq:heat_flux_BC}
\end{equation}
allow us to model a prescribed heat flux.
In fact there is a specific form of heat flux that is equivalent to the constant Dirichlet condition $f(t)=T_0$ in Eq.~(\ref{BCs}), that is~\cite{Bouciguez_Thermal_Engineering_2006}
\begin{equation}
    \frac{\partial T}{\partial x} (0,t) = -  \frac{q_0}{k_L \sqrt{t}}  \, .
    \label{eq:special_heat_flux_BC}
\end{equation}
Hence, given a relation between $q_0$ and $T_0$, the analytic solution Eq.~(\ref{stefanSol}) is applicable also in this case,
as we illustrate numerically in the upcoming Sect.

The implementation of Dirichlet boundary conditions was described in Sect.~\ref{sec:RW_BC}.
Here we sketch an implementation of the Neumann boundary condition~(\ref{eq:heat_flux_BC}).

At the first time step, we seed the temperature for the fixed boundary with the order of unity, i.e. $T(x_0,t_0) \simeq 1 $.  
Using a forward differentiation approximation
\begin{equation}
    \frac{\partial T}{\partial x} (0,t) \approx   \frac{ T(x_1,t_j) - T(x_0,t_j)}{ \Delta x}  \, ,
\end{equation}
we in the consecutive time steps ($j>0$) update the temperature at the fixed boundary according to 
\begin{equation}
 T(x_0,t_j) = \text{round} \left( -n \Delta x h(t_j) +  T(x_1,t_j) \right) \, ,
\end{equation}
where \textit{round} rounds a number to the nearest integer.

\section{Numerical results for Stefan problems}\label{results}
\begin{figure}[p]
    \includegraphics[width = 60mm]{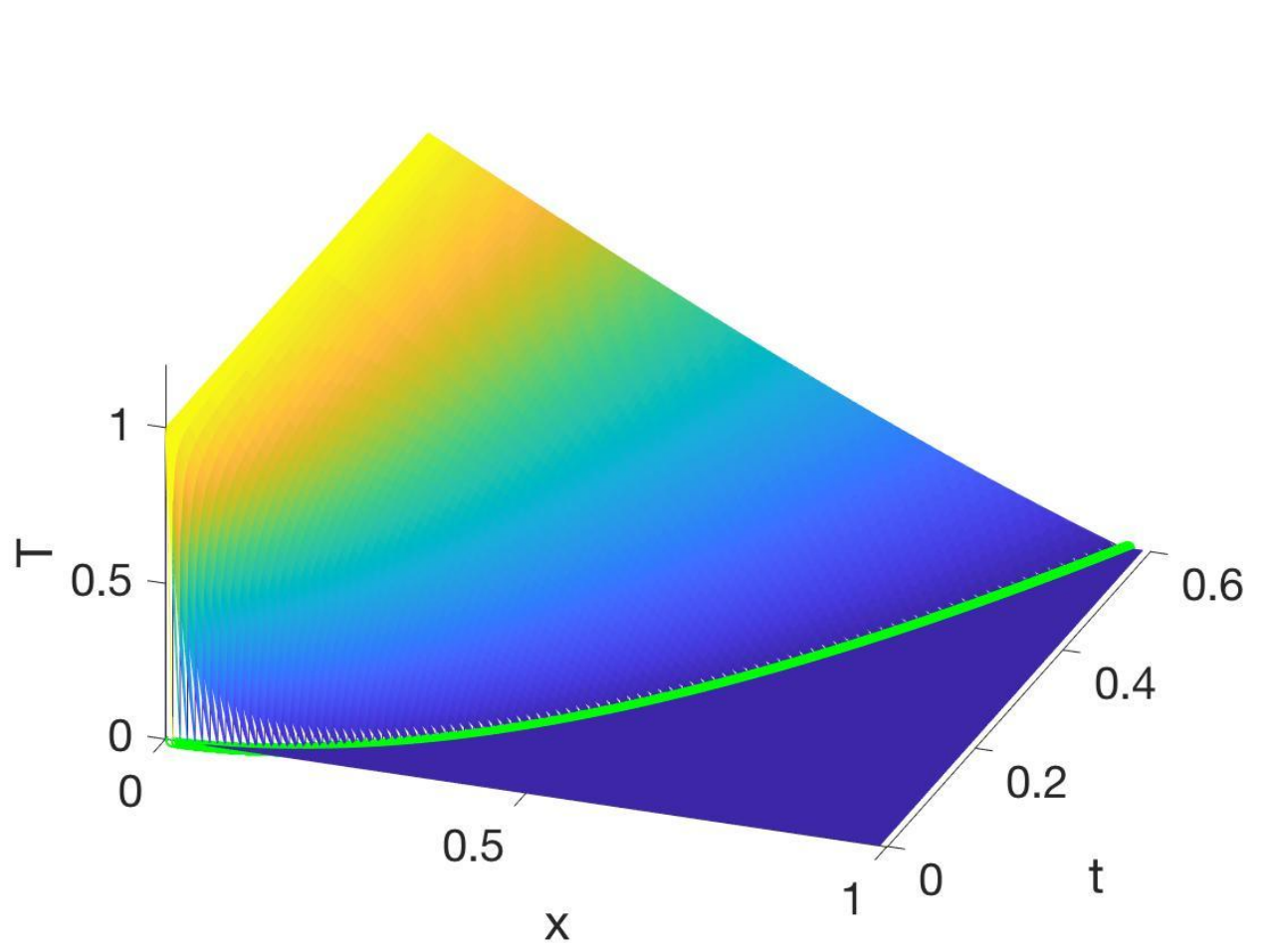}~\includegraphics[width = 60mm]{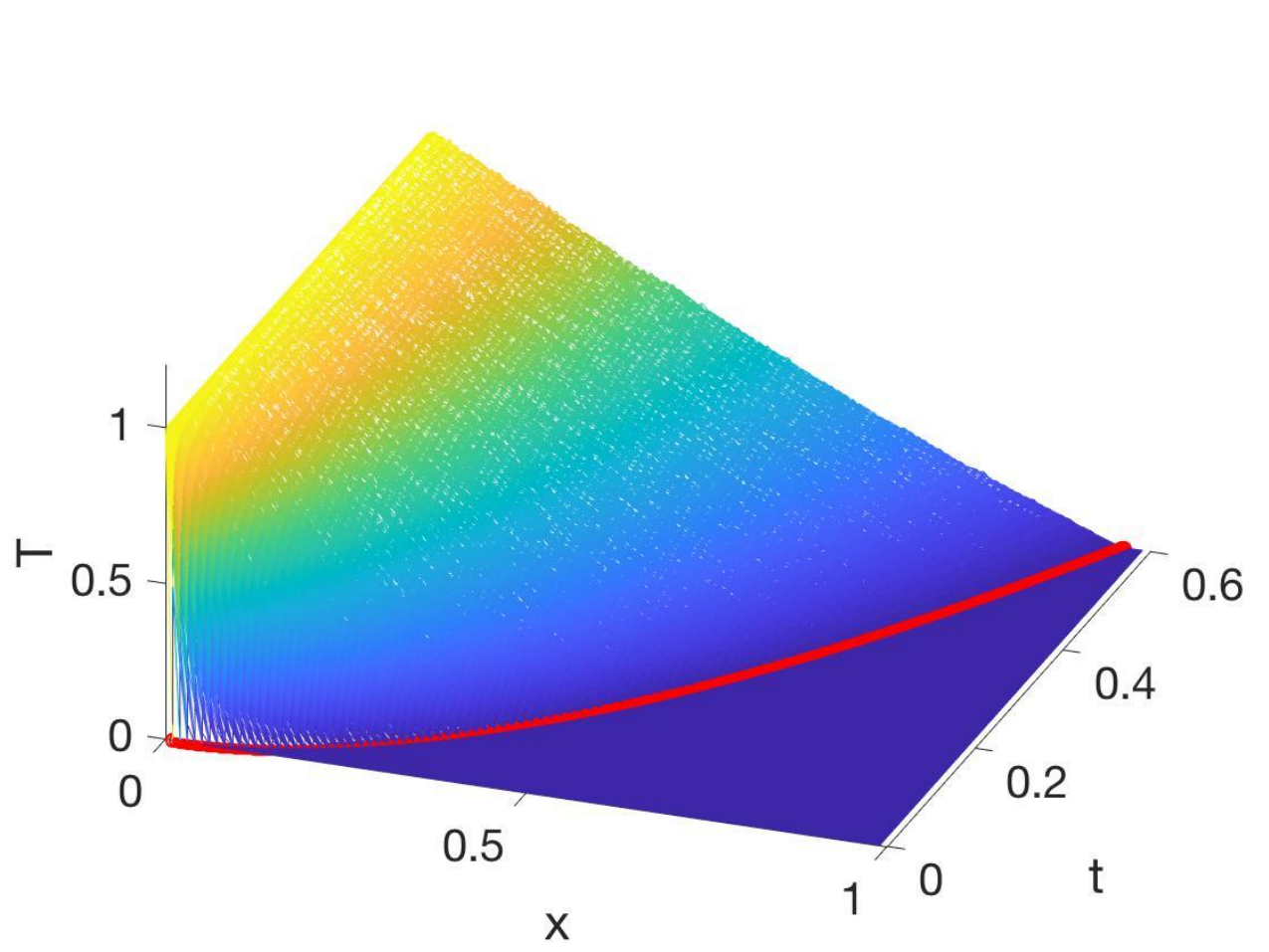}
    \caption{Solutions for Stefan problem for $t \in [0,0.6]$ with boundary condition $f(t)=1$. 
    (Left) Analytic solution from Eq.~(\ref{stefanSol}). 
    (Right) RWM solution with $n=10^4$ and $\Delta x=0.01$.}
    \label{fig_const_b}
\end{figure}
\begin{figure}[p]
\sidecaption[t]
    \includegraphics[width = 75mm]{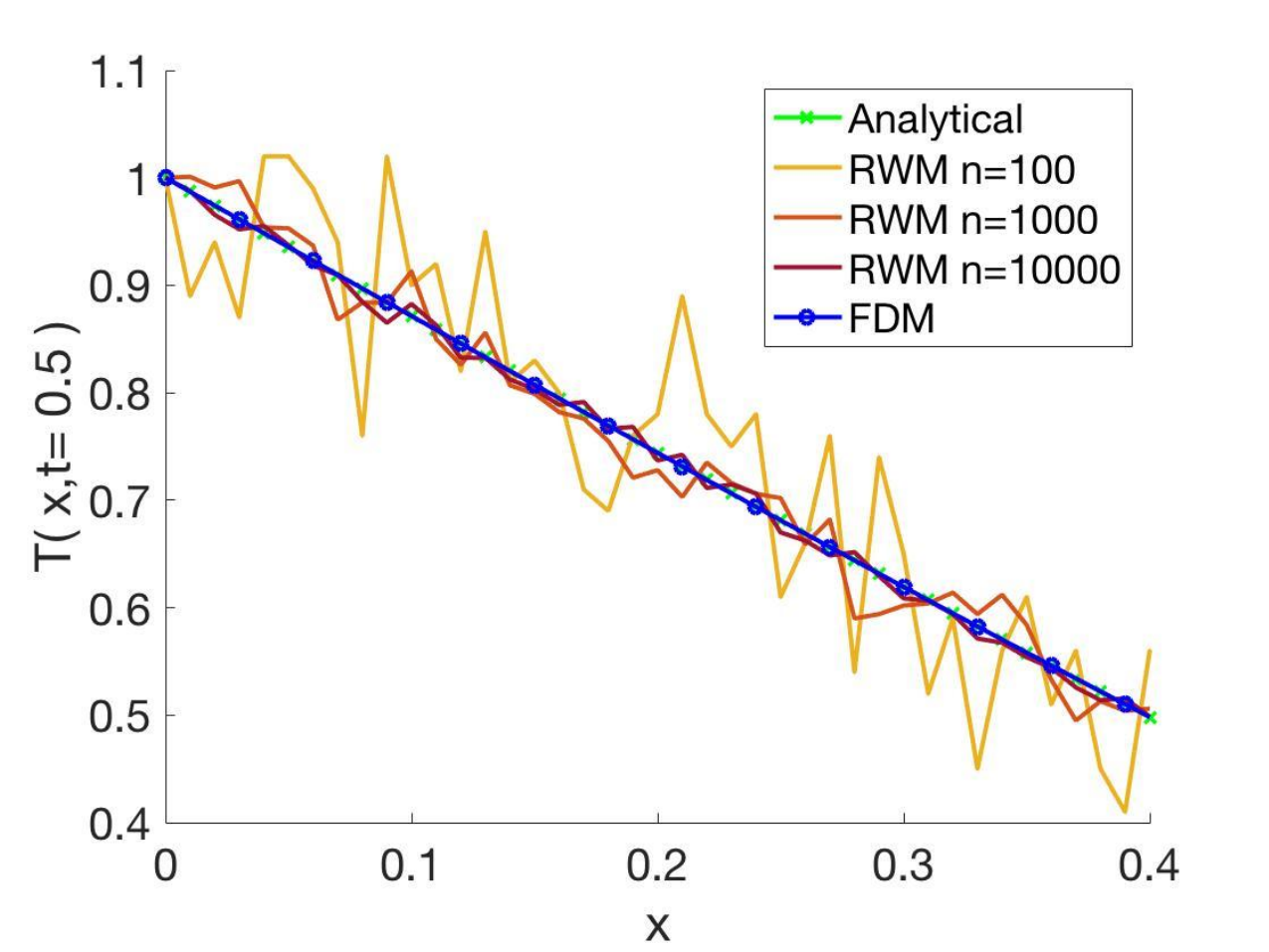}
    \caption{Analytic solution $T(x,0.5)$, $x \in [0,0.4]$, of Eq.~(\ref{stefanSol}) (green). 
    RWM solutions with $\Delta x=0.01$ and different values of $n$, see inset legend. 
    Numerical results from a finite difference method (FDM)~\cite{umea} (blue).}
    \label{const_left}
\end{figure}
\begin{figure}[p]
\sidecaption[t]
    \includegraphics[width = 75mm]{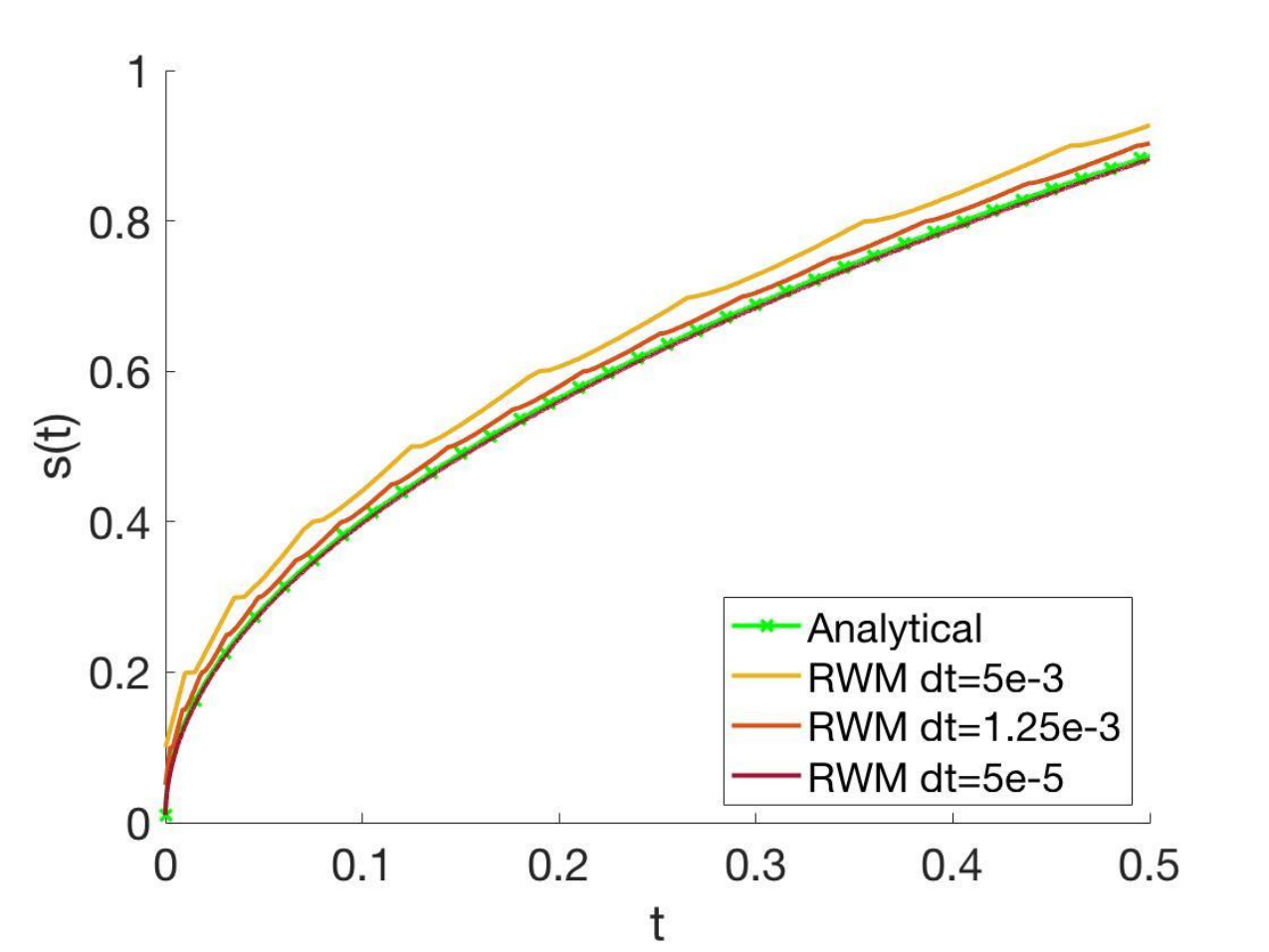}
    \caption{ 
    Analytic solution $s(t)$ of Eq.~(\ref{stefanSol}) (green). 
    Numerical solution of $s(t)$ from RWM for $t \in [0,0.5]$ and $n=10^4$, with different step lengths $\Delta t$, see inset legend.
    }
    \label{const_right}
\end{figure}

\subsection{Stefan problem with constant boundary condition $f(t)=T_0$}

In Eq.~(\ref{stefanSol}) we presented the analytic solution for the Stefan problem Eqs.~(\ref{heat1})-(\ref{ste1}) when $f(t)=T_0$. Fig. \ref{fig_const_b} shows the temperature distributions $T(x,t)$ for the analytic result and the RWM solution with $T_0=1^\circ$C, for $\alpha=1$ and $\beta=1$. The green respectively the red curves denotes the solid-liquid interface. 
In Fig.~\ref{const_left} we compare different values of $n$ for the RWM in the cross-section $t=0.5$. 
In the Fig.~\ref{const_right} we compare different sizes of the step length $\Delta t$ in a plot of the moving boundary $s(t).$

\begin{figure}[t]
    \includegraphics[width = 60mm]{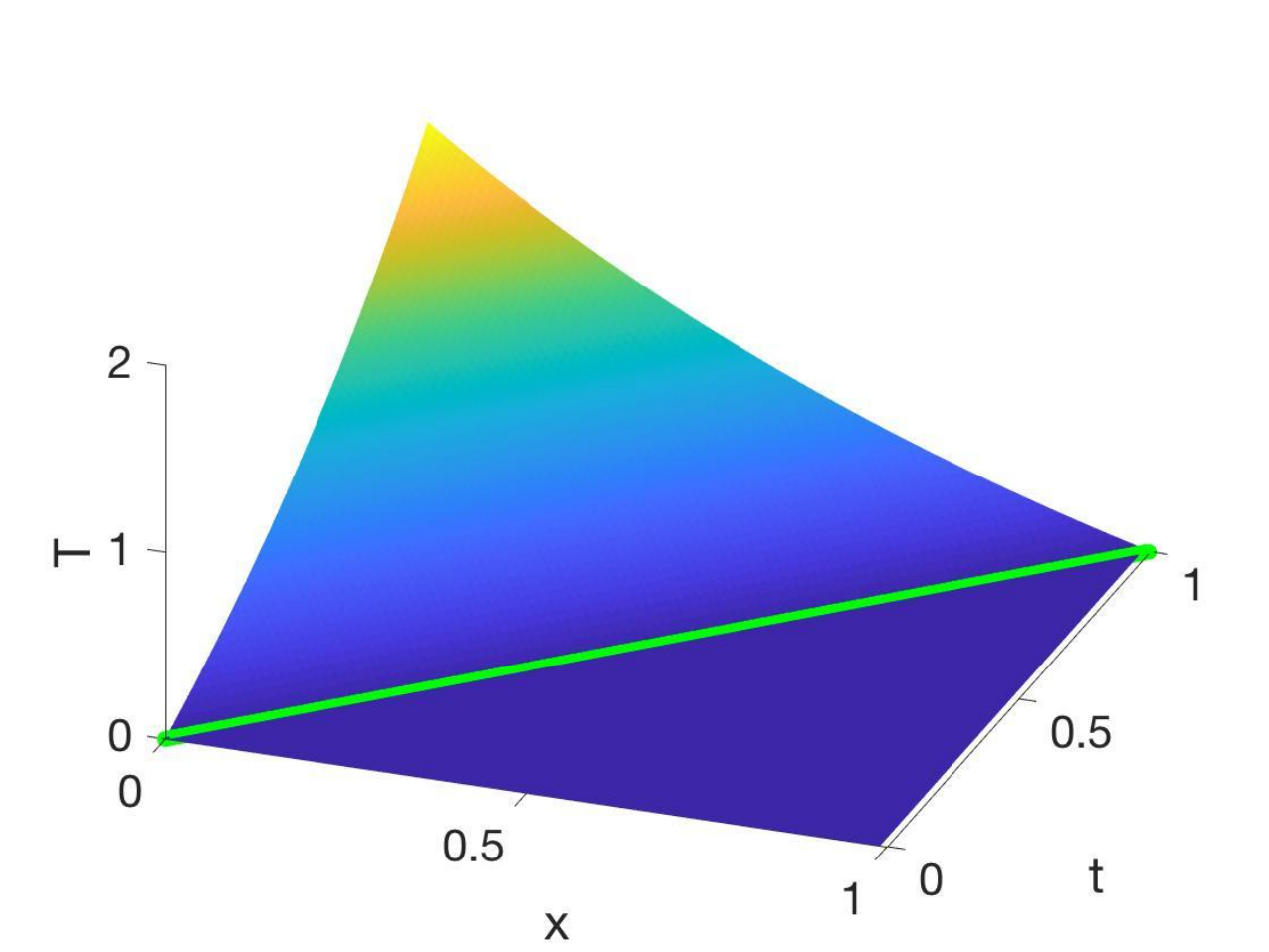}~\includegraphics[width = 60mm]{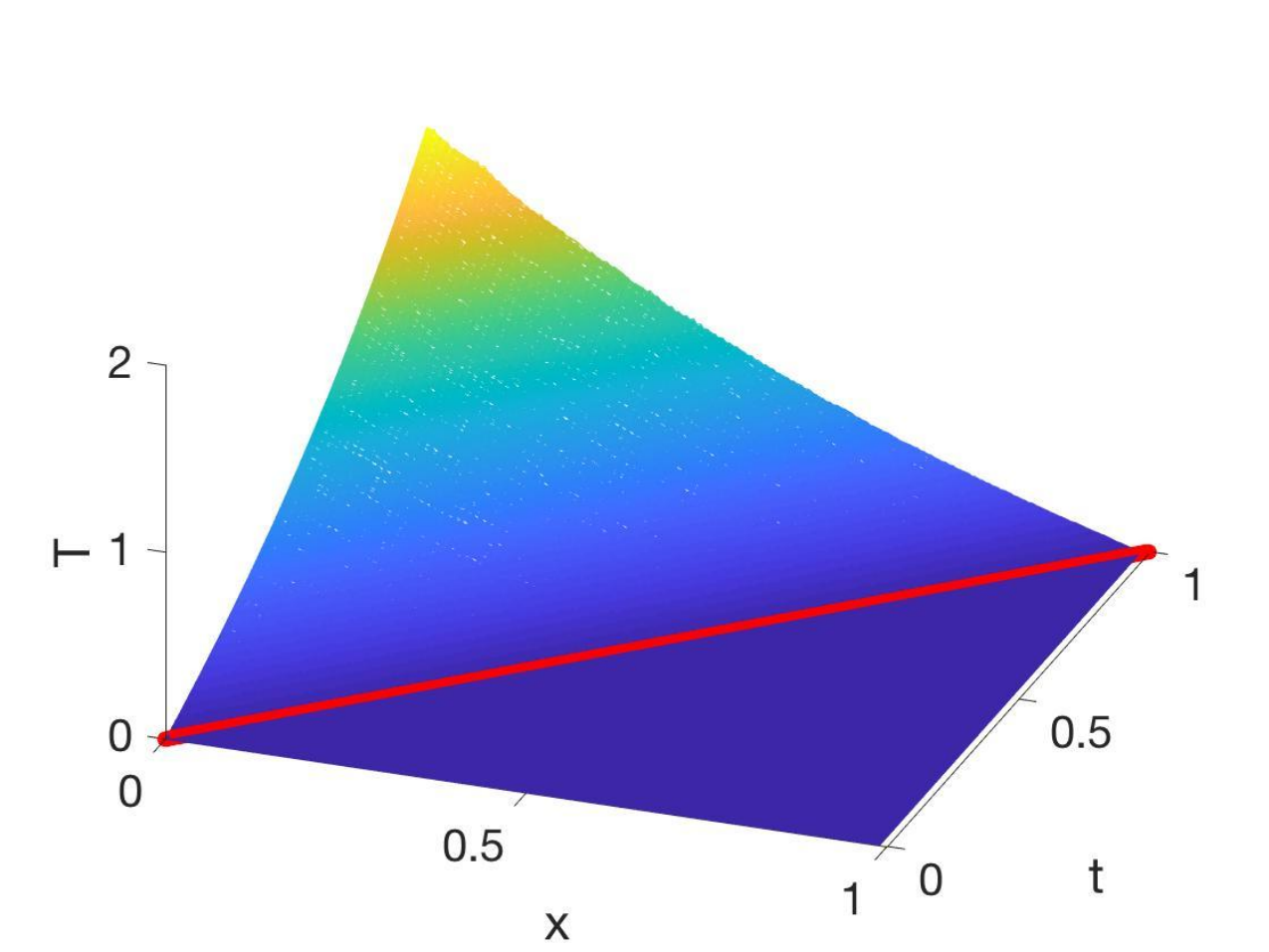}
    \caption{Solutions for the Stefan problem with the boundary condition $f(t)=e^t-1$ for $t \in [0,1]$.
    (Left) Analytic solution Eq.~(\ref{special}). 
    (Right) RWM solution with $n=10^4$ and $\Delta x=0.01$.}
    \label{fig_exp}
\end{figure}
\begin{figure}[t]
\sidecaption[t]
          \includegraphics[width=75mm]{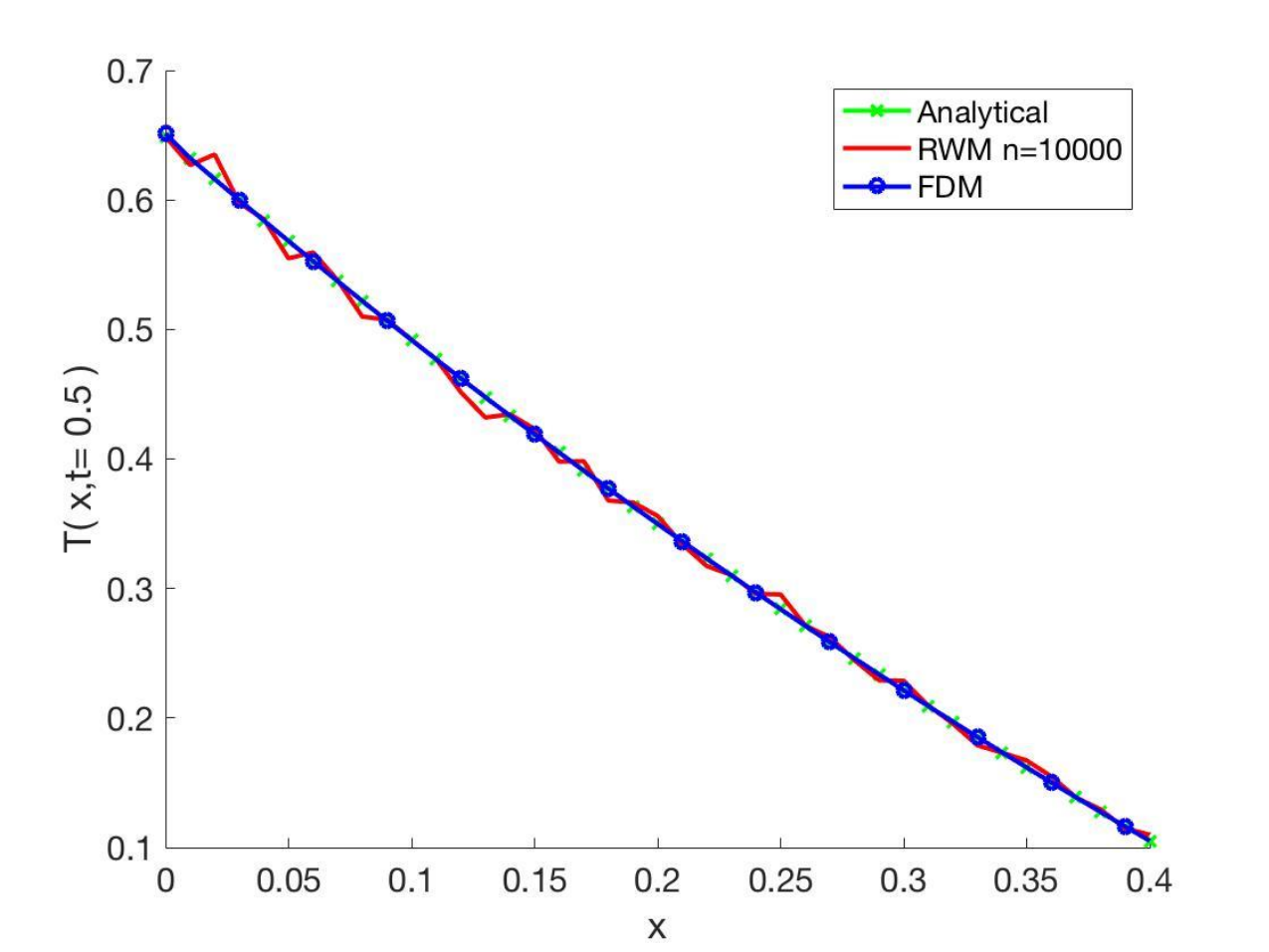}
    \caption{Analytic solution $T(x,0.5)$, $x \in [0,0.4]$, of Eq.~(\ref{special}) (green). 
    RWM solutions with $n=10^4$ and $\Delta x=0.01$ (red). 
    Numerical results from a finite difference method (FDM)~\cite{umea} (blue).}
    \label{fig_exp_b}
\end{figure}

\subsection{Stefan problem with a special boundary condition $f(t)=e^t-1$}

In Eq.~(\ref{special}) we presented the analytical solution for the Stefan problem Eqs.~(\ref{heat1})-(\ref{ste1}) in the special case when $f(t)=e^t-1$. 
Fig.~\ref{fig_exp} shows the temperature distributions $T(x,t)$ for the analytical result and the RWM solution for $\alpha=1$ and $\beta=1$. 
The green respective the red curves shows the solid-liquid interface. 
In Fig.~\ref{fig_exp_b} we see a comparison between the analytic result, the RWM and FDM in the cross-section $t=0.5$.

\subsection{Stefan problem with a special heat flux boundary condition $h(t)=-q_0/\sqrt{t}$}

\begin{figure}[t]
    \includegraphics[width = 60mm]{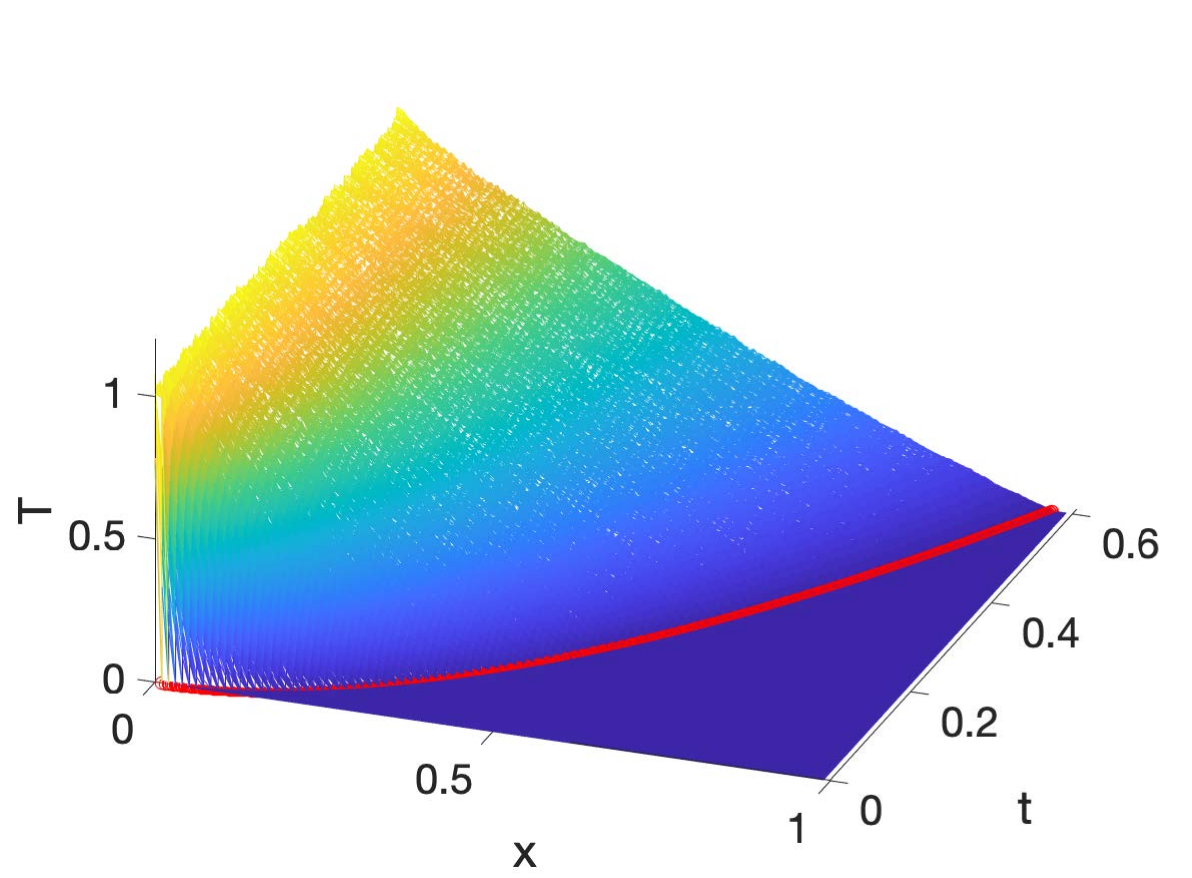}~\includegraphics[width = 50mm]{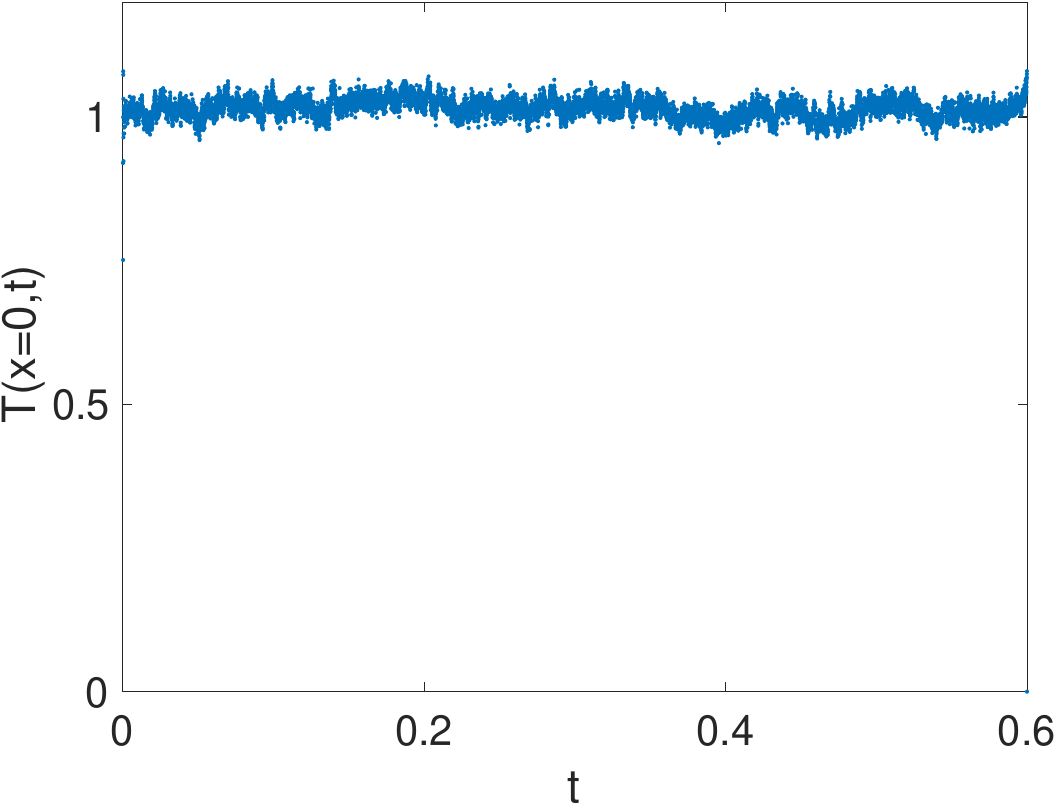}
  \caption{Solutions of the Stefan problem with the special heat flux boundary condition $h(t)= - q_0/ \sqrt{t}$ for $t \in [0,0.6]$.
    (Left) RWM solution for $q_0=0.9108$ with $n=10^4$ and $\Delta x=0.01$, to be compared with Fig.~\ref{fig_const_b}. 
    (Right) Cross section $x=0$ of the RWM solution.}
    \label{fig_heat_flux}
\end{figure}

We now estimate what value of $q_0$ that is required in order for the temperature to be $T(0,t)=T_0$.
The total heat entering during the time $t$ is
 \begin{equation}
   Q = S \frac{q_0}{k_L} \int_0^t \frac{d \tau}{\sqrt{\tau}} = \frac{2 S q_0}{k_L} \sqrt{t} \, .
       \label{eq:Q_heat_flux}
\end{equation}

From Fig.~\ref{const_left} we obtain the approximation $T(x,t=0.5) \approx 1 - 1.25 x$ for the constant temperature case.
Hence, during the time interval $t \in [0, 0.5]$, the solid phase have received the heat $Q_S = \rho V l   \approx  \rho  S (s(0.5) - s(0)) l = 0.8 \rho l$,
and the liquid phase have received the heat $Q_L = \rho  V  c \Delta T  \approx \rho S (s(0.5) - s(0))  c (T(0,0.5) + T(0.4,0.5))/2    = 0.6 \rho c $.
With $\rho  = S = l = c = k_L = 1$ ($\alpha  = \beta = 1$), we have the total heat $Q = Q_S + Q_L = 1.4$.
Solving for $q_0$ from Eq.~(\ref{eq:Q_heat_flux}), we obtain the estimation $q_0 \approx 0.99$.
If one instead calculates $\partial T(0,t) / \partial x$ from the analytic solution Eq.~(\ref{stefanSol}), 
one obtains $q_0=1/(\sqrt{\pi} \text{erf}(\lambda)) = 0.9108$.
Numerically we find that $q_0 \approx 0.9108$ gives a constant temperature $T(x=0,t) \approx  1$ for $\alpha=1$ and $\beta=1$, see Fig.~\ref{fig_heat_flux},
which is in agreement with~\cite{Bouciguez_Thermal_Engineering_2006}.

\subsection{Stefan problem with oscillating boundary condition}

  \begin{figure}[t] 
      \includegraphics[width=60mm]{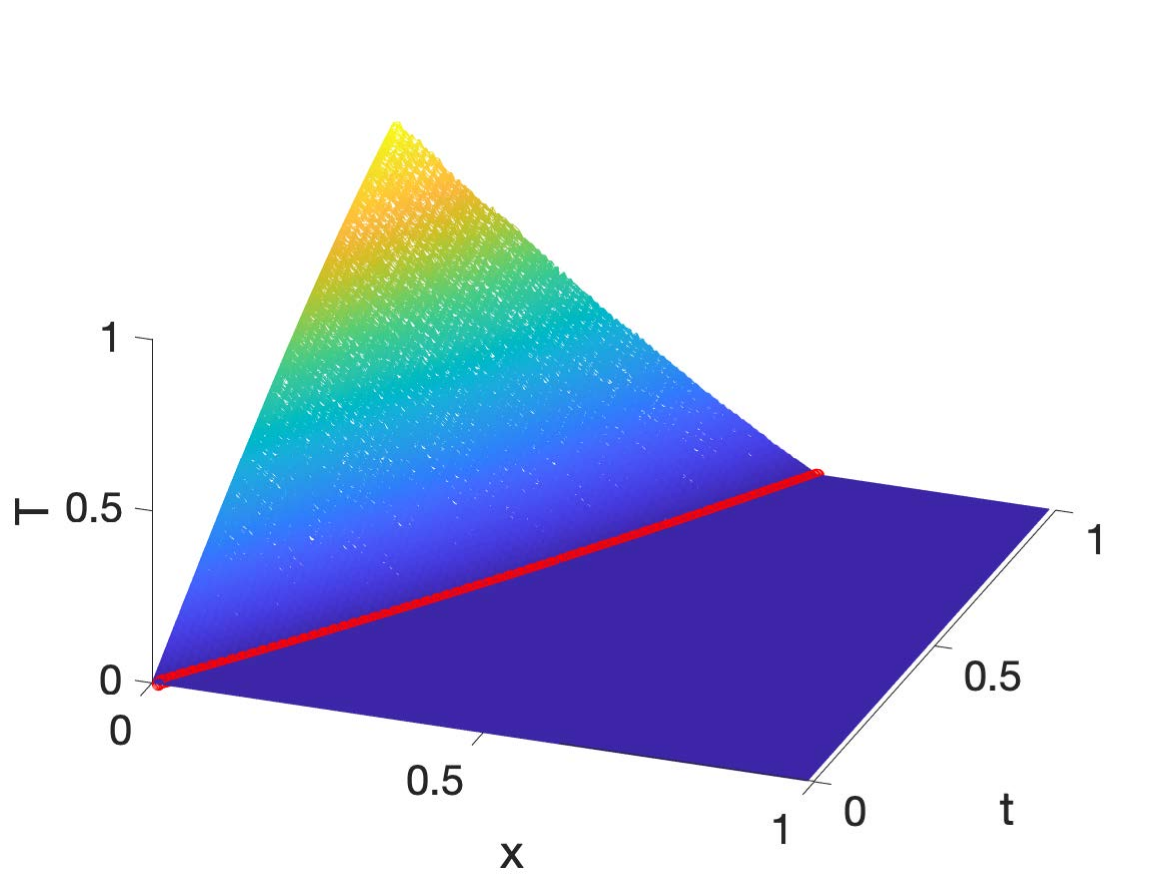}~\includegraphics[width=60mm]{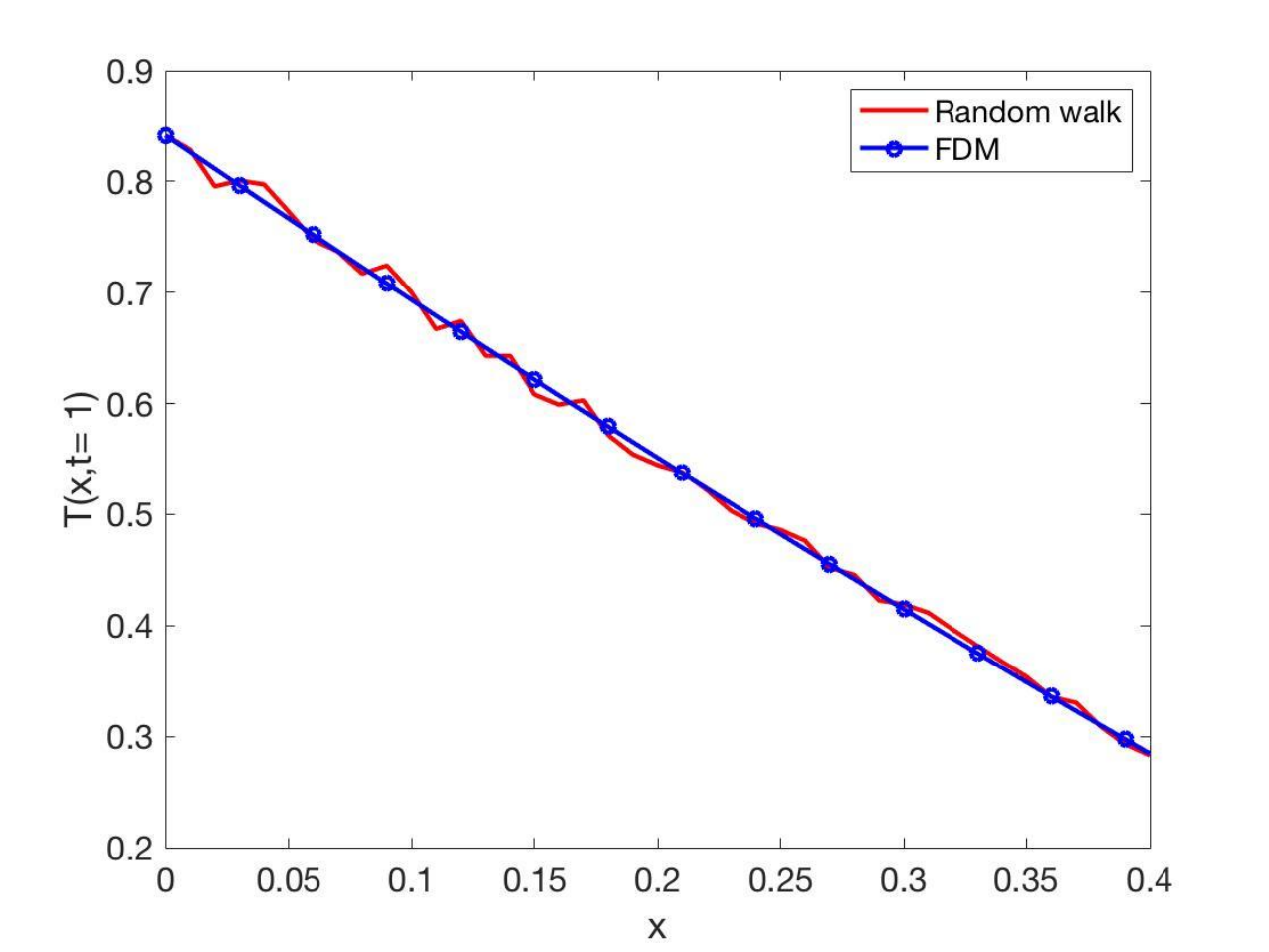}
      \caption{(Left) RWM solution for $f(t)=\sin{( t)}$ with $n=10^4$ and $\Delta x=0.01$ on the interval $t \in [0,1]$.
      (Right) Comparison between RWM and FDM for $f(t)=\sin(t)$ at $t=1$ on the interval $x \in [0,0.4]$.}
      \label{fig_sin1}
  \end{figure}
  
In the introduction we proposed to model a general time dependent fixed boundary condition $T(0,t)=f(t)$ with the RWM.
Due to limitations in the existing code for the finite difference method (FDM)~\cite{umea}, we are presently restricted to consider $f(t)=\sin (t)$ at the boundary when comparing the two numerical methods.
The RWM solution for the Stefan problem Eqs.~(\ref{heat1})-(\ref{ste1}) yields the temperature distribution $T(x,t)$ 
as seen in the left part of Fig.~\ref{fig_sin1}.
The RWM is compared to the FDM for the cross-section $t=1$ in the right part of Fig.~\ref{fig_sin1}.  
Here we have set $\alpha =1$ and $\beta=2$.  
\begin{figure}[h]
\sidecaption[t]
      \includegraphics[width=75mm]{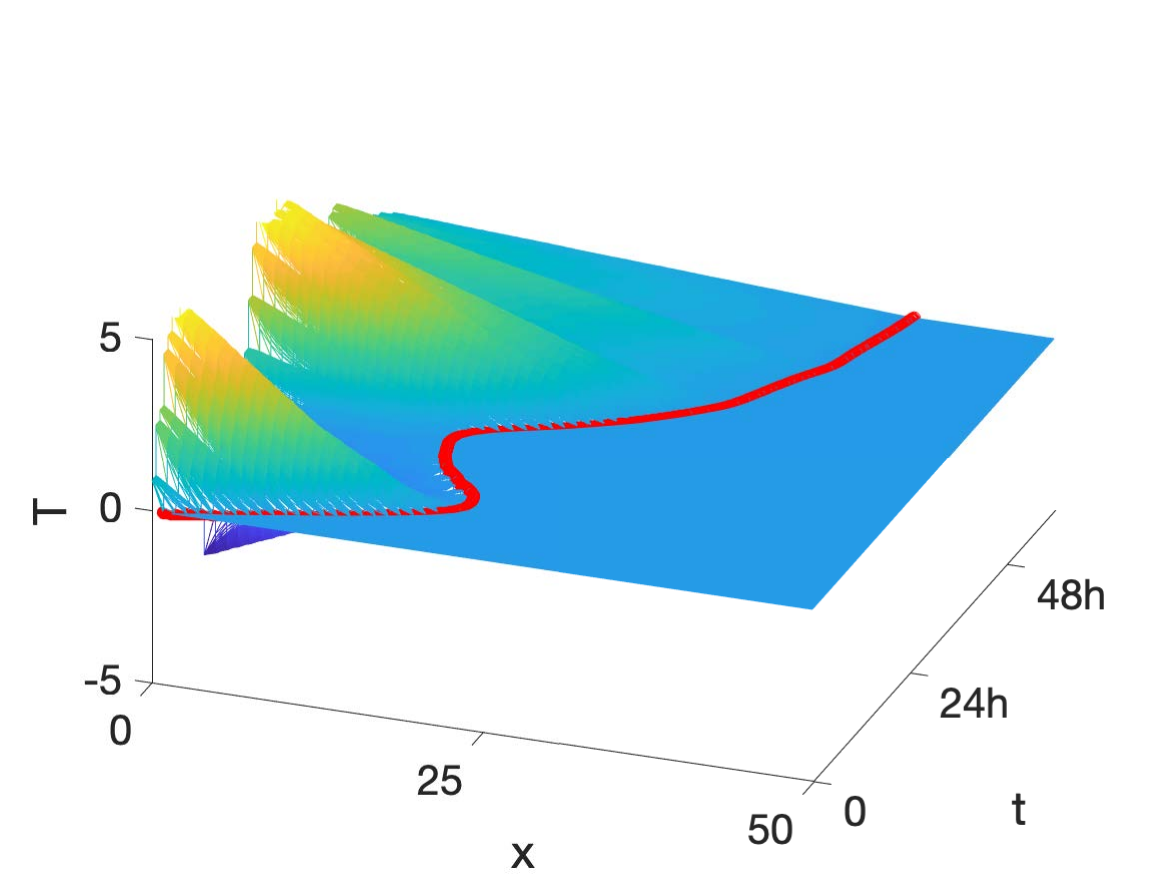}
      \caption{RWM model where $f(t)$ is set to the observed day temperatures at \"{O}rebro airport 1-3 March 2019, $x$ is in mm and $t \in [0, 62 \, \text{h}]$.}
      \label{fig_day1}
\end{figure}
\subsection{Stefan problem with boundary condition according to daytime temperature variations}

To finally apply our RWM model with an arbitrary time dependent temperature at the fixed boundary in a simulation of melting ice, 
we set the physical constants for water to $\alpha \approx 0.1429 \text{ mm}^2/\text{s}$ and $\beta \approx 79.9 \text{ K}$~\cite{lockby}. 
We model the melting of ice according to the daytime temperature variations and therefore we set $f(t)$ at the fixed boundary to the observed air temperatures from \"{O}rebro airport 1-3 March 2019~\cite{smhi}. 
Assuming the observed air temperature at the fixed surface is a simplification that does not take the temperature gradient between air and ice/water, or heat transport by convection or radiation, into account. 
Nevertheless, Fig.~\ref{fig_day1} gives a qualitative view of the dynamics of the melting ice, and we see for example that it is freezing again during the first night, although the present one-phase implementation with negative temperatures in the liquid is quantitatively unrealistic.

\section{Discussion}

From the numerical results of the previous Sect., we can see qualitatively from Figs.~\ref{fig_const_b}-\ref{fig_heat_flux} that the RWM solution to Stefan problems converges to the analytical as $\Delta x \rightarrow 0$. 
An oscillatory boundary condition was successfully evaluated against a finite difference method in Fig.~\ref{fig_sin1}.
Finally, an arbitrary time dependent function for the fixed boundary was used to model the melting of ice with realistic temperature data in Fig.~\ref{fig_day1}.

There are a few simplifications in our model for the Stefan problem that can be improved in a more detailed study. 
Among the physical simplifications, we have mentioned our assumption that we use the same density for water and ice, $\rho_L=\rho_S$ which is not the real case. 
We may also want to consider a temperature distribution in the solid phase, $T_{ice}(x,t) \neq 0$, which leads to a two-phase Stefan problem with a system of PDE:s. 
Some cases of two-phase problems also have analytic solutions, see e.g.~\cite{furzeland}.

There are several applications for the Stefan problems in different fields of engineering. 
By looking at the original purpose of Stefan's article in 1891, which was to model the arctic ices, this is highly relevant today due to the demand of better climate models. 
According to Hunke \textit{et al.}~\cite{hunke}, Stefan's one dimensional thermodynamical model is still in use for global climate models, although the complete thermodynamical sea-ice models are of course more complex. 
Hence thermodynamical sea-ice models may be a subject for future work with the RWM approach. 

Other areas where a solid-liquid interface is moving is in 3D-printing, freezing of food, solidifying of building components.
Also, in lithium-ion batteries, the diffusion of lithium ions in the battery is separated into two phases, one where lithium ions are evenly distributed, and one where they are not present. 
To be able to compute the properties of batteries in a better way, such as life-time and capacity, 
one can estimate the movement of the interface between these two phases
as a Stefan problem~\cite{hoffman}.

\section{Conclusions}
In accordance with our opening objective, we have successfully used a stochastic method to calculate numerical solutions with arbitrary accuracy to the Stefan problem with general time-dependent boundary condition at the fixed boundary.
In comparison with the finite difference method, 
our experience is that the RWM is easier to implement and 
more flexible in terms of switching between different boundary conditions.
This further motivates the use of stochastic methods in more complex applied problems in higher dimensions~\cite{OgrenJMR2019}.

\begin{acknowledgement}
We thank the students Andreas Lockby, Daniel Stoor, and Emil Gestsson for fruitful discussions about the Stefan problem. 
We are also grateful to Tobias Jonsson for sharing the finite difference code, used here for comparisons with the stochastic method.
Finally we thank Daniel Edstr\"{o}m and Bair Budaev for proofreading.
\end{acknowledgement}

%\newpage

\section*{Appendix}
\addcontentsline{toc}{section}{Appendix}

\begin{programcode}{Example of a \textsc{Matlab} code that can reproduce Figs.~\ref{fig_const_b}-\ref{fig_sin1}} %% that can Program Code}
%If you want to emphasize complete paragraphs of texts in an \verb|Program Code| we recommend to
%use

%\verb|\begin{programcode}{Program Code}|

%\verb|\begin{verbatim}...\end{verbatim}|

%\verb|\end{programcode}|

\begin{verbatim}
% RWM_Stefan.m (can be downloaded from the arXiv:2006.04939 [math.AP] Ancillary files)
clear all; close all
% PARAMETERS:
alpha=1 % K/(rho*c); % Thermal diffusivity.
beta=1 % l/c; % Parameter with unit [K]. 
L=1 % Length of domain
t_max=0.5 % Maximum time
T_0=1; % [degree C] Temperature for constant temperature BC.
% Parameter for the constant heat flux BC.
q_0 = 0.9108 % = 1/(sqrt(pi)*erf(lambda)).
n=1e2; % Number of walkers.
dx=0.01; dt=dx^2/(2*alpha); % Steplengths in x and t
ds=dx/(n*beta); % Increment for s(t) when absorbing a walker.
% Number of points in the space and time.
N_x=ceil(L/dx); N_t=ceil(t_max/dt); 
% Matrix representing T(x,t), initially set to 0 degree C.
T=zeros(N_x,N_t); 
s_vector=zeros(1,N_t); % Vector representing s(t).
j_t=1; j_s=1; % Indices for time and the position of s(t).
s=dx; % initial value for s(t) /approx 0.
% Loop for all time steps as long as s(t) < L.
while j_t < N_t && j_s < N_x 
% Examples of boundary conditions (BC) for the fixed boundary.
T(1,j_t)=n*T_0; % Constant Dirichlet BC.
% T(1,j_t)=n*(exp(j_t*dt)-1); % Exponential BC.
% T(1,j_t)=n*sin(j_t*dt); % Oscillating BC.
% % Heat flux
% if j_t==1 % First timestep.
%    T(1,1)=n*1; % Seed temperature of order unity.
% else % Consecutive timesteps.
%    T(1,j_t)=round( (n*dx*q_0/(j_t*dt)^(0.5)+T(2,j_t)) ); 
% end % if
s_vector(j_t)=s;
for j_x = 1:N_x
    if T(j_x,j_t) < 0 % If T is below 0 degree C (unrealistic one-phase model).
        sign=-1;
    else
        sign=1;
    end
    for k=1:sign*T(j_x,j_t) % Move all walkers at (j_x,j_t).
        p=2*round(rand)-1; % =+-1, with P(+1)=P(-1)=1/2.
        % A walker move if it has not reached the boundaries.
        if j_x+p > 1 && j_x+p <= j_s && j_x <= N_x 
            T(j_x+p,j_t+1) = T(j_x+p,j_t+1) + 1*sign;
        elseif j_x+p == j_s+1 % Otherwise s(t) moves ds. 
            s=s+ds*sign; % Update s(t).
            j_s=floor(s/dx); % New index for s(t).
        end
    end % k
end % j_x
% For calculating the heat flux.
q_0_vector(j_t)=(T(1,j_t)-T(2,j_t))/dx; 
j_t=j_t+1; % Update the time index.
end % while
T=T/n; % Dividing by the number of walkers.
% Plot the temperature distribution and s(t).
figure; hold on 
[x_matrix,t_matrix]=meshgrid(0:dx:(N_x-1)*dx,0:dt:(N_t-1)*dt);
mesh(x_matrix',t_matrix',T)
t_vector=0:dt:(j_t-2)*dt;
plot3(s_vector(1:j_t-1),t_vector, 0*t_vector,'ro')
xlabel('x'); ylabel('t'); zlabel('T')
view([20 40]); set(gca,'FontSize',20)
% Plot s(t) from the RWM and the analytical solution.
figure; hold on 
% For constant T or special heat flux BC, the analytical  
% solution requires solution of the transcendental equation. 
lambda=trans_eq(beta,T_0); 
s_ana_vector = dx+2*lambda*sqrt(alpha*t_vector);
plot(t_vector, s_vector(1:j_t-1),'r.')
% The solution for constant BC.
plot(t_vector, s_ana_vector,'g--','Linewidth',2) 
% The solution for exponential BC (beta=1).
%plot(t_vector, t_vector,'g--','Linewidth',2) 
xlabel('t'); ylabel('s(t)'); set(gca,'FontSize',20)
% Solving for lambda from the transcendental equation 
% with Newton-Raphson method:
function [x0] = trans_eq(beta,T_0)
  f=@(x) sqrt(pi)*beta*x*exp(x^2)*erf(x)-T_0;
  fprim=@(x) beta*(sqrt(pi)*exp(x^2)*erf(x)*(2*x^2+1)+2*x^2);
  tol=1e-6; % Tolerance.
  x0=1; % initial guess.
  while abs(f(x0)) > tol
     x0 = x0-f(x0)/fprim(x0);
  end
end

\end{verbatim}

\end{programcode}

\end{document}